\documentclass{birkmult}

\usepackage{amsmath,amsfonts}
\usepackage{amsthm}
\usepackage{enumerate}

\newtheorem{theorem}{Theorem}[section]
\newtheorem{proposition}[theorem]{Proposition}

\theoremstyle{definition}
\newtheorem{definition}[theorem]{Definition}
\newtheorem{example}[theorem]{Example}
\newtheorem{condition}[theorem]{Condition}

\newtheorem{remark}[theorem]{Remark}

\newtheorem*{remark*}{Remark}

\numberwithin{equation}{section}

\newcommand{\nC}{\mathbb C}
\newcommand{\nR}{\mathbb R}

\newcommand{\cD}{{\mathcal D}}
\newcommand{\cF}{{\mathcal F}}
\newcommand{\cK}{{\mathcal K}}
\newcommand{\cL}{{\mathcal L}}

\newcommand{\kip}{[\,\cdot\, , \cdot\,]}

\newcommand{\co}[1]{\overline{#1}}
\newcommand{\V}[1]{\mathbf{#1}}
\newcommand{\M}[1]{\mathsf{#1}}

\DeclareMathOperator{\sgn}{sgn} 
\DeclareMathOperator{\dom}{\cD} \DeclareMathOperator{\fdom}{\cF}

\begin{document}
\title[Riesz Bases of Root Vectors, I]%
{Riesz Bases of Root Vectors of \\
Indefinite Sturm-Liouville Problems \\
with  Eigenparameter  Dependent \\
Boundary Conditions, I}

\author{Paul Binding}
\address{
Department of Mathematics and Statistics, University of Calgary,
Calgary, Alberta, Canada, T2N 1N4}
 \email{binding@ucalgary.ca}
\author{Branko \'{C}urgus}
\address{Department of Mathematics,
Western Washington University, Bellingham, WA 98225, USA}
 \email{curgus@cc.wwu.edu}



\date{\today}

\begin{abstract}
We consider a regular indefinite Sturm-Liouville problem with two
self-adjoint boundary conditions, one being affinely dependent on
the eigenparameter.   We give sufficient conditions under which a
basis of each root subspace for this Sturm-Liouville problem can
be selected so that the union of all these bases constitutes a
Riesz basis of a corresponding weighted Hilbert space.
\end{abstract}

\maketitle

\section{Introduction} \label{s2}

We consider a regular indefinite Sturm-Liouville boundary
eigenvalue problem of the form
\begin{equation} \label{sl1}
-(p\,f')'+q\,f = \lambda \, r\, f \ \ \ \ \text{on} \ \ \ \
[-1,1].
\end{equation}
The coefficients $1/p, q, r$ in \eqref{sl1} are assumed to be real
and integrable over $[-1,1]$, $p(x) > 0$, and $x\, r(x) > 0$ for
almost all $x \in [-1,1]$.  We impose two boundary conditions on
\eqref{sl1} (only one of which is $\lambda$-dependent):
\begin{equation} \label{bc1}
\M{L}\, \V{b}(f) = \M{0} , \ \ \ \ \ \ \  \M{M} \V{b}(f) = \lambda\,
\M{N} \V{b}(f).
\end{equation}
where $\M{L},\M{M}$ and $\M{N}$ are $1 \times 4$ nonzero (row)
matrices and the boundary mapping $\V{b}$ is defined for all $f$
in the domain of \eqref{sl1} by
\begin{equation*}
\V{b}(f) = \begin{bmatrix} f(-1) & f(1) & (pf')(-1) &
(pf')(1)\end{bmatrix}^T .
\end{equation*}

We shall utilize an operator theoretic framework developed in
\cite{BC}.  Under Condition~\ref{cbc} below, a self-adjoint
operator $A$ in the Krein space $L_{2,r}(-1,1)\oplus \nC_{\Delta}$
can be associated with the eigenvalue problem \eqref{sl1},
\eqref{bc1}. Here $\Delta$ is a nonzero real number which is
determined by $\M{M}$ and $\M{N}$ -- see Section~\ref{sop} for
details.  We remark that the topology of this Krein space is that
of the corresponding Hilbert space $L_{2,|r|}(-1,1)\oplus
\nC_{|\Delta|}$. (In the rest of the paper we abbreviate
$L_{2,r}(-1,1)$ to $L_{2,r}$ and $L_{2,|r|}(-1,1)$ to
$L_{2,|r|}$.)  Our main goal in this paper is to provide
sufficient conditions on the coefficients in \eqref{sl1},
\eqref{bc1} under which there is a Riesz basis of the above
Hilbert space consisting of the union of bases for all the root
subspaces of the above operator $A$. This will be referred to for
the remainder of this section as the {\em Riesz-basis property of}
$A$.

Completeness and expansion theorems with a stronger topology, but
in a smaller space corresponding to the form domain of the
operator $A$, have been considered by many authors -- see
\cite{BC} (and the references there) and \cite{LS}. Although the
topology of the Krein space $L_{2,r}\oplus \nC_{\Delta}$ is weaker
than the topology of the form domain, which in our case is a
Pontryagin space, the expansion question turns out to be much more
challenging mathematically.

Indeed, even for the case when the boundary conditions are
$\lambda$-independent this problem is nontrivial. In our notation,
this case corresponds to $\M{L}$ being a nonsingular $2 \times 4$
matrix, with the second equation in \eqref{bc1} suppressed.  The
Riesz-basis property of the operator corresponding to $A$, now
defined in $L_{2,r}$, has been discussed by several authors, e.g.,
in \cite{Be,CL,F1,Py,V}.  The first general sufficient condition for
this was given by Beals \cite{Be}, who required the weight function
$r$ to behave like a power of the independent variable $x$ in an
open neighborhood of the turning point $x=0$, although his method
does allow more general weight functions. Refinements of Beals's
method in \cite{F1} and \cite{V} show that a ``one-sided" condition
on $r$ (i.e., in only a half-neighborhood of $x=0$ ) is enough to
guarantee the Riesz-basis property.  That some extra condition on
$r$ is indeed necessary follows from \cite{V} where Volkmer showed
that weight functions $r$ exist for which the corresponding
Sturm-Liouville problem \eqref{sl1}, under the conditions used here,
does not have the Riesz-basis property. Explicit examples of such
weight functions were given in \cite{AP} and \cite{F2}. Recently,
Parfyonov \cite{Pa} has given an explicit necessary and sufficient
condition for the Riesz basis property in the case $p = 1, q = 0$
with odd weight function $r$. Here, and in most of the above
references, Dirichlet boundary conditions were imposed.

General self-adjoint (perhaps non-separated, but still
$\lambda$-independent) boundary conditions were treated by
\'{C}urgus and Langer \cite{CL}.  They showed that if the {\em
essential} boundary conditions, i.e., those not including
derivatives, were separated, then a Beals-type condition in a
neighborhood of $x=0$ was sufficient for the Riesz-basis property.
But if some of the non-separated boundary conditions were
essential then \cite{CL} established the Riesz-basis property only
by imposing extra restrictions on the weight function in
(half-)neighborhoods of {\em both} endpoints of the interval
$[-1,1]$.  Again, some extra restriction is necessary, since in
\cite{BC1} we gave an explicit example of \eqref{sl1} under the
conditions used here, satisfying a Beals-type condition at $x=0$,
but without the Riesz-basis property. Of course at least one (in
fact one, in this antiperiodic case) boundary condition was
essential and non-separated. In some sense, then, the boundary
$\pm 1$ behaves as a turning point under such boundary conditions.

In summary, the Riesz-basis property is quite subtle, and depends
significantly on the nature of the boundary conditions even when
they are independent of $\lambda$. In this paper and its sequel,
we shall examine the analogous situation for the cases of one and
two $\lambda$-dependent boundary conditions, where the
possibilities for the ($\lambda$-dependent) boundary conditions
are much greater.  As in the $\lambda$-independent case, a
condition on the weight function is needed near the turning point
$x=0$ to ensure the Riesz-basis property of $A$. We shall develop
such a condition (which is implied by the ones discussed above) in
Section~\ref{sc0}.  Depending on the nature of the boundary
conditions \eqref{bc1}, we may also need a condition near the
boundary, and this is discussed in Section~\ref{scm11}. It should
be remarked that for the case of exactly one $\lambda$-dependent
boundary condition treated here we need only one such condition,
near either $x=-1$ or $x=1$, and this can be viewed as a
``one-sided" condition at $\pm 1$.  In the case of two
$\lambda$-dependent boundary conditions we shall also need a
condition involving both boundary points $x=-1$ and $x=1$.

It turns out that all the above conditions have a common core.
This is not immediately obvious, since there are differences
between the ``turning points" $0$ and $\pm 1$. For example, when
the boundary conditions are separated, the values of $f$ and $f'$
are equal at $0$ but are independent at $-1$ and $1$.  The common
core, which will also be needed in Part II, involves the notion of
smoothly connected half-neighborhoods, and this is defined and
studied in Section~\ref{ssc}.

In order to apply the above conditions, we use a criterion in
Theorem~\ref{W}, equivalent to the Riesz-basis property of $A$,
involving a positive homeomorphism of $L_{2,r}\oplus \nC_{\Delta}$
with the form domain of $A$ as an invariant subspace.  This,
together with certain mollification arguments, is used for our
main results, which are detailed in Section~\ref{sk=0}. To
paraphrase these, we recall that a $\lambda$-independent boundary
condition is {\em essential} if it does not include derivatives.
Similarly, a $\lambda$-dependent boundary condition will be called
{\em essential} if it does not include derivatives in the
$\lambda$-terms.

In Theorem~\ref{tk0sebc} we discuss situations when a condition on
$r$ near $x=0$ suffices for the Riesz-basis property of $A$. For
example this holds when the first ($\lambda$-independent) boundary
condition in \eqref{bc1} is either non-essential, or essential and
separated, and the second ($\lambda$-dependent) one is
non-essential. If the latter condition is essential instead, then
the same result holds if a sign condition is also satisfied, and
this includes a result of Fleige \cite{F3}, which is the only
reference we know where the Riesz-basis property of $A$ has been
studied for $\lambda$-dependent boundary conditions.

In Theorems~\ref{td0k04} and \ref{tn=k2} we consider those cases
of \eqref{bc1} which are not covered by Theorem~\ref{tk0sebc}.
Then we require a condition near just one of the boundary points
$\pm 1$, not both as in \cite{CL}. The choice of the boundary
point is arbitrary in Theorem~\ref{td0k04} which deals with the
case when the boundary conditions in \eqref{bc1} are,
respectively, essential non-separated and non-essential. In
Theorem~\ref{tn=k2}, however, this choice is not arbitrary but
depends on the sign of the number $\Delta$ used in defining the
inner product on $L_{2,r}\oplus \nC_{\Delta}$.

\section{Operators associated with the eigenvalue problem}
\label{sop}

The maximal operator $S_{\max}$ in $L_{2,r}(-1,1)=L_{2,r}$
associated with \eqref{sl1} is defined by
\[
S_{\max} : f \mapsto \ell(f) := \frac{1}{r} \bigl(-(pf')' + qf
\bigr), \ \ \ f \in \cD(S_{\max}),
\]
where
\[
\dom(S_{\max}) =  \cD_{\max} = \bigl\{ f \in L_{2,r} : f, pf' \in
AC[0,1], \ \ell(f) \in L_{2,r} \bigr\}.
\]
We define the boundary mapping $\V{b}$ by
\begin{equation*}
\V{b}(f) := \begin{bmatrix} f(-1) & f(1) & (pf')(-1) &
(pf')(1)\end{bmatrix}^T, \ \ \  f \in \dom(S_{\max}),
\end{equation*}
and the concomitant matrix $\M{Q}$ corresponding to $\V{b}$ by
\[
\M{Q} := \text{{\Large $i$}} \begin{bmatrix} 0  &  0 &  - 1 & 0  \\
0 & 0 & 0 & 1  \\ 1 & 0 & 0 & 0  \\ 0 & -1  & 0 & 0
\end{bmatrix}.
\]
We notice that $\M{Q} = \M{Q}^{-1}$. Integrating by parts we
easily calculate that
\[
\int_{-1}^1 S_{\max}f \, \co{g}\, r  - \int_{-1}^1 f\,
S_{\max} \co{g} \, r  =  i \, \V{b}(g)^* \M{Q} \V{b}(f),
\ \ \  f, g \in \dom(S_{\max}).
\]

Throughout, we shall impose the following nondegeneracy and
self-adjointness conditions on the boundary data.

\begin{condition} \label{cbc}
The row vectors $\M{L}, \M{M}$ and $\M{N}$ in \eqref{bc1}
satisfy:
\begin{enumerate}[(1)]
\item
the $3 \times 4$ matrix $\begin{bmatrix} \M{L} \\ \M{M} \\
\M{N}\end{bmatrix}$ has rank $3$,
\item
 $\M{L} \M{Q} \M{L}^* = \M{M
Q M}^* = \M{N Q N}^* = \M{L Q M}^* = \M{L Q N}^* = 0$,
\item \label{ibp4}
 $i\, \M{M} \M{Q}^{-1} \M{N}^*$ is a nonzero real number
 and we define
\begin{equation} \label{Del}
\M{\Delta} = -  \frac{i}{\M{M} \M{Q}^{-1} \M{N}^*}.
\end{equation}
\end{enumerate}
\end{condition}

Clearly the boundary value problem \eqref{sl1}-\eqref{bc1} will
not change if row reduction is applied to the coefficient matrix
\begin{equation} \label{cm}
\begin{bmatrix} \M{L} & \M{0} \\ \M{M} & \M{N}  \end{bmatrix}.
\end{equation}
In what follows we will assume that the $2 \times 8$ matrix in
\eqref{cm} is row reduced to row echelon form (starting the
reduction at the bottom right corner). After the row reduction,
we write the row vectors $\M{L}$ and $\M{N}$ as
\begin{equation} \label{Le}
\M{L} =  \begin{bmatrix} \M{L}_{e}  &  \M{L}_n \end{bmatrix} \, ,
\ \ \ \ \ \ \ \M{N} =  \begin{bmatrix} \M{N}_{e}  & \M{N}_n
\end{bmatrix}.
\end{equation}
If either of the $1 \times 2$ matrices $\M{L}_n, \M{N}_n$ is
nonzero, the corresponding boundary condition is called
``non-essential".  In any case these matrices do not appear in the
representation of the form domain of $A$, discussed below, but
they will play an important role in our conditions for Riesz bases
in Section~\ref{sk=0}. The $1 \times 2$ matrices $\M{L}_{e}$ and
$\M{N}_{e}$ represent the ``essential" boundary conditions if the
non-essential parts $\M{L}_{n}$ and $\M{N}_{n}$ are zero matrices.

Next we define a Krein space operator associated with the problem
\eqref{sl1}-\eqref{bc1}.  We consider the linear space $L_{2,|r|}
\oplus \nC$, equipped with the inner product
\[
\left[ \begin{pmatrix} f \\ z \end{pmatrix} , \begin{pmatrix} g
\\ w \end{pmatrix} \right] := \int_{-1}^{1} f \co{g} r \, +
\, \co{w} \Delta z\, , \ \ \  f, g \in L_{2,|r|}, \ z, w \in \nC.
\]
Then $\bigl(L_{2,|r|} \oplus \nC,\kip\bigr)$ is a Krein
space, which we denote by $L_{2,r} \oplus \nC_{\Delta}$.
A fundamental symmetry on this Krein space is given by
\begin{equation} \label{fsy}
J := \begin{bmatrix} J_0  & 0  \\  0 & \sgn\M{\Delta}
\end{bmatrix},
\end{equation}
where $\sgn\M{\Delta} \in \{-1,1\}$ and $J_0:L_{2,r} \to L_{2,r}$
is defined by
\[
 (J_0f)(x) :=  f(x)\, \sgn(r(x)), \ \ \ x \in [-1,1].
\]
Then $[J\,\cdot\,,\,\cdot\,]$ is a positive definite inner product
which turns $L_{2,r} \oplus \nC_{\Delta}$ into a Hilbert space
$L_{2,|r|} \oplus \nC_{|\Delta|}$.  The topology of $L_{2,r}
\oplus \nC_{\Delta}$ is defined to be that of $L_{2,|r|} \oplus
\nC_{|\Delta|}$, and a Riesz basis of $L_{2,r} \oplus
\nC_{\Delta}$ is defined as a homeomorphic image of an orthonormal
basis of $L_{2,|r|} \oplus \nC_{|\Delta|}$.

We define the operator $A$ in the Krein space $L_{2,r} \oplus
\nC_{\Delta}$ on the domain
\begin{equation} \label{domwT}
\dom(A) =  \left\{ \begin{bmatrix} f \\ z
\end{bmatrix}  \in
 \begin{matrix} L_{2,r} \\ \oplus \\ \nC_{\Delta} \end{matrix}
  : f \in \dom\bigl(S_{\max}\bigr),
\ \M{L}\, \V{b}(f)=0, \ \ z = \M{N} \V{b}(f)  \right\}
\end{equation}
by
\begin{equation} \label{defwT}
A  \begin{bmatrix} f \\ \M{N} \V{b}(f)
\end{bmatrix} :=
\begin{bmatrix} S_{\max}f \\ \M{M} \V{b}(f) \end{bmatrix}, \ \
\ \ \ \ \ \ \begin{bmatrix} f \\ \M{N} \V{b}(f)
\end{bmatrix} \in \dom(A).
 \end{equation}

Using \cite[Theorems 3.3 and 4.1]{BC} we see that this operator is
self-adjoint in $L_{2,r} \oplus \nC_{\Delta}$ and in particular:
\begin{enumerate}[(i)]
 \item
$A$ is quasi-uniformly positive \cite{CN} (and therefore
definitizable) in $L_{2,r} \oplus \nC_{\Delta}$.
\item
$A$ has a discrete spectrum.
\item
The root subspaces corresponding to real distinct eigenvalues of
$A$ are mutually orthogonal in the Krein space
$L_{2,r}\oplus\nC_{\Delta}$.
\item
All but finitely many eigenvalues of $A$ are semisimple and real.

\end{enumerate}

For further properties of $A$, we refer the reader to
\cite[Theorem 3.3]{BC}. From (i), (ii) and the characterization of
the regularity of the critical point infinity for definitizable
operators in Krein spaces given in \cite[Theorem 3.2]{C}, we then
obtain the following, which is our central tool.

\begin{theorem} \label{W}
Let $\fdom(A)$ denote the form domain of $A$. There exists a basis
for each root subspace of $A$, so that the union of all these
bases is a Riesz basis of $L_{2,|r|} \oplus \nC_{|\Delta|}$ if and
only if there exists a bounded, boundedly invertible, positive
operator $W$ in $L_{2,r} \oplus \nC_{\Delta}$ such that $W
\fdom(A) \subset \fdom(A)$.
\end{theorem}

In order to apply this result, we need to characterize the form
domain $\fdom(A)$.  To this end, let $\cF_{\max}$ be the set of
all functions $f$ in $L_{2,r}$, absolutely continuous on $[-1,1]$,
such that $\int_{-1}^1 p\, |f^{\prime}|^2 < +\infty$. On
$\cF_{\max}$ we define the {\em essential boundary mapping}
$\V{b}_e:\cF_{\max} \rightarrow \nC^{2}$ by
\begin{equation*} 
 \V{b}_e(f) := \begin{bmatrix} f(-1) & f(1) \end{bmatrix}^{T}, \
\ \ f \in  \cF_{\max}.
\end{equation*}
Clearly $\V{b}_e$ is surjective.

By \cite[Theorem 4.2]{BC}, there are four possible cases for the
form domain $\fdom(A)$ of $A$: If $\M{L}_n \neq \M{0}$ and
$\M{N}_n \neq \M{0}$, then
\begin{align}
\fdom(A) & = \left\{ \begin{bmatrix} f \\ z \end{bmatrix} \in
 \begin{matrix} L_{2,r} \\ \oplus \\ \nC_{\Delta} \end{matrix}
  : f \in \cF_{\max}, \  z \in \nC
\right\}.\label{fd1}
 \intertext{
If $\M{L}_n = \M{0}$ and $\M{N}_n \neq \M{0}$, then }
 \fdom(A)
& = \left\{ \begin{bmatrix} f \\ z \end{bmatrix} \in
\begin{matrix} L_{2,r} \\ \oplus \\ \nC_{\Delta} \end{matrix} : f
\in \cF_{\max}, \M{L}_e \V{b}_e(f) = 0, \ z \in \nC
\right\}.\label{fd2}
 \intertext{
If $\M{L}_n \ne \M{0}$ and $\M{N}_n = \M{0}$, then  }
 \fdom(A) & = \left\{ \begin{bmatrix} f \\ \M{N}_e \V{b}_e(f)
\end{bmatrix} \in
\begin{matrix} L_{2,r} \\ \oplus \\ \nC_{\Delta} \end{matrix}
: f \in \cF_{\max} \right\}.\label{fd3}
 \intertext{
If $\M{L}_n = \M{0}$ and $\M{N}_n = \M{0}$, then  }
 \fdom(A) & = \left\{ \begin{bmatrix} f \\ \M{N}_e \V{b}_e(f)
\end{bmatrix} \in
 \begin{matrix} L_{2,r} \\ \oplus \\ \nC_{\Delta} \end{matrix}
  : f \in \cF_{\max}, \M{L}_e \V{b}_e(f) =
0 \right\}.\label{fd4}
\end{align}

To construct an operator $W$ as in Theorem~\ref{W} we need to
impose conditions on the coefficients $p$ and $r$ in \eqref{sl1}.
In all cases we need Condition~\ref{c0} in a neighborhood of $0$,
and in some cases we also need one of two Conditions, \ref{cat-1}
or \ref{cat1}, on $r$ in neighborhoods of $-1$ or $1$. These will
be discussed in Sections~\ref{sc0} and \ref{scm11} respectively.

\section{Smooth connection and associated operator}
\label{ssc}

To prepare the ground for the Conditions mentioned above (and in
Part II), we develop the concept of smoothly connected
half-neighbourhoods. A closed interval of non-zero length is said to
be a {\em left half-neighborhood} of its right endpoint and a
{\em right half-neighborhood} of its left endpoint.

Let $\imath$ be a closed subinterval of $[-1,1]$.  By
$\cF_{\max}(\imath)$  we denote the set of all functions $f$ in
$L_{2,r}(\imath)$ which are absolutely continuous on $\imath$ and
such that $\int_{\imath} p\, |f^{\prime}|^2 < +\infty$.  Note that
$\cF_{\max}[-1,1]$ is the space $\cF_{\max}$ defined below Theorem
\ref{W}. In the next definition affine function $\alpha$ means
$\alpha(t) = a+\alpha't$ where $a,\alpha', t \in \nR$.

\begin{definition} \label{dscab}
Let $p$ and $r$ be the coefficients in \eqref{sl1}. Let $a, b \in
[-1,1]$ and let $h_a$ and $h_b$, respectively,  be
half-neighborhoods of $a$ and $b$ which are contained in $[-1,1]$.
We say that the ordered pair $(h_a, h_b)$ is {\em smoothly
connected} if there exist
\begin{enumerate}[(a)\ \ ]
 \item
positive real numbers $\epsilon$ and $\tau$,
 \item
non-constant affine functions $\alpha:[0,\epsilon] \rightarrow
h_a$ and $\beta:[0,\epsilon] \rightarrow h_b$,
 \item
non-negative real functions $\rho$ and $\varpi$ defined on
$[0,\epsilon]$,
\end{enumerate}
such that
\begin{enumerate}[(i)\ \ ]
\item \label{idscabi}
$\alpha(0) = a$ \ and \ $\beta(0) = b$,
 \item \label{idscabii}
$p\circ\alpha$ and $p\circ\beta$ are locally integrable on the
interval $(0,\epsilon]$,
 \item \label{idscabiii}
$\rho\circ \alpha^{-1} \in
\cF_{\max}\bigl(\alpha([0,\epsilon])\bigr)$,
\item \label{idscabiiia}
$1/\tau < \varpi < \tau$ a.e. on $[0,\epsilon]$,
 \item \label{idscabiv}
$\displaystyle \rho(t)
  = \frac{\,\bigl|r\bigl(\beta(t)\bigr)\bigr|}%
  {\bigl|r\bigl(\alpha(t)\bigr)\bigr|}$,
  \ \  and \ \
  $\displaystyle   \varpi(t)  =
\frac{\,p\bigl(\beta(t)\bigr)}{p\bigl(\alpha(t)\bigr)}$, \ \ \
 for  \ \ \ $\displaystyle  t \in (0,\epsilon]$.
\end{enumerate}
The numbers $\alpha', \beta'$, (the slopes of $\alpha$, $\beta$,
respectively) and $\rho(0)$ are called the {\em parameters} of the
smooth connection.
\end{definition}

\begin{remark}
Since the function $\alpha$ in Definition~\ref{dscab} is affine,
the condition $\rho\, \circ\, \alpha^{-1} \in
\cF_{\max}\bigl(\alpha([0,\epsilon])\bigr)$ in (\ref{idscabiii})
is equivalent to
\begin{equation} \label{eqdscabiii}
\rho \in AC[0,\epsilon] \ \ \ \ \text{and} \ \ \ \
 \int_0^{\epsilon} |\rho'(t)|^2 p(\alpha(t))dt < + \infty.
\end{equation}
Under the assumption that $1/\tau < \varpi < \tau$ a.e.
on $[0,\epsilon]$, it follows that property \eqref{eqdscabiii}
is equivalent to
 \[
\rho \in AC[0,\epsilon] \ \ \ \ \text{and} \ \ \ \
 \int_0^{\epsilon} |\rho'(t)|^2 p(\beta(t))dt < + \infty.
 \]
\end{remark}

To illustrate Definition~\ref{dscab}, we make the following
\begin{definition}
Let $\nu$ and $a$ be real numbers and let $h_a$ be a
half-neighborhood of $a$.  Let $g$ be a function defined on $h_a$.
Then $g$ is called {\em of order} $\nu$ {\em on} $h_a$ if there
exists $g_1 \in C^1(h_a)$ such that
\begin{equation*}
 g(x) = (x-a)^{\nu}g_1(x) \ \ \ \ \text{and} \ \ \ \ g_1(x) \neq
 0, \ \ \ x \in h_a.
\end{equation*}
\end{definition}

\begin{example}
Let $a, b \in [-1,1]$. Let $h_a$ and $h_b$, respectively,  be
half-neighborhoods of $a$ and $b$ contained in $[-1,1]$. Assume
that the coefficient $r$ in \eqref{sl1} is  of order $\nu$ on both
half-neighborhoods $h_a$ and $h_b$.  Assume also that the
functions $p$ and $1/p$ are bounded on $h_a$ and $h_b$ (or,
alternatively, that $p$ is  of order $\mu$ on both
half-neighborhoods $h_a$ and $h_b$.)  Then lengthy, but
straightforward, reasoning shows that the half-neighborhoods $h_a$
and $h_b$ are smoothly connected. Moreover the parameters of the
smooth connection are nonzero numbers.
\end{example}

\begin{remark} \label{rd0}
Throughout the paper we use the following convention: A product of
functions is defined to have value $0$ whenever one of its terms
has value zero, even if some other terms are not defined.
\end{remark}

\begin{theorem} \label{tgenS}
Let $\imath$ and $\jmath$ be closed intervals, $\imath, \jmath \in
\bigl\{[-1,0],[0,1] \bigr\}$. Let $a$ be an endpoint of $\,\imath$
and let $b$ be an endpoint of $\,\jmath$.  Denote by $a_1$ and
$b_1$, respectively, the remaining endpoints.  Assume that the
half-neighborhoods $\imath$ of $a$ and $\jmath$ of $b$ are
smoothly connected with parameters $\alpha', \beta'$ and
$\rho(0)$.  Then there exists an operator
 $$
S: L_{2,|r|}(\imath)\rightarrow L_{2,|r|}(\jmath)
 $$
such that:
\begin{enumerate}[{\rm ($S$-1)\ }]
\item \label{itS1}
$S \in \cL\bigl(L_{2,|r|}(\imath),L_{2,|r|}(\jmath)\bigr), \ S^*\in
\cL\bigl(L_{2,|r|}(\jmath),L_{2,|r|}(\imath)\bigr)$.
\item \label{itS2}
$(Sf)(x) = 0, \ |x - b_1| \leq \frac{1}{2}$ for all $f \in
L_{2,|r|}(\imath)$, \ and \\
$(S^*g)(x) = 0, \ |x - a_1| \leq \frac{1}{2}$ for all $g \in
L_{2,|r|}(\jmath)$.
\item \label{itS4}
$S\cF_{\max}(\imath) \subset \cF_{\max}(\jmath)$, \ \
$S^*\cF_{\max}(\jmath) \subset \cF_{\max}(\imath)$.
\item \label{itS3}
For all $f \in \cF_{\max}(\imath)$ and all $g \in
\cF_{\max}(\jmath)$ we have
\[
\lim\limits_{\substack{y\rightarrow b \\ y \in \jmath\,}} \
(Sf)(y) \, = \, |\alpha'| \,
 \lim\limits_{\substack{x\rightarrow a\\ x \in
\imath}} \, f(x), \ \ \  \lim\limits_{\substack{x\rightarrow a\\
x \in \imath}} \ (S^*g)(x) \ = \ |\beta'|\rho(0) \,
\lim\limits_{\substack{y\rightarrow b \\ y \in \jmath\,}} \,
g(y).
\]
\end{enumerate}
\end{theorem}
\begin{proof}
Let $\epsilon > 0$ be the real number and $\alpha$ and
$\beta$ the affine functions introduced in
Definition~\ref{dscab}.  Thus $\alpha(0) = a$ and
$\beta(0) = b$. It is no loss of generality to assume
that each of the intervals
$\alpha\bigl([0,\epsilon]\bigr)$ and
$\beta\bigl([0,\epsilon]\bigr)$ has a length $ < 1/2$.
Let $\alpha_1:[0,1] \to \imath$ and $\beta_1:[0,1] \to
\jmath$ be strictly monotonic and continuously
differentiable bijections such that $\alpha_1(x) =
\alpha(x)$ and $\beta_1(x) = \beta(x)$ for all $x \in
[0,\epsilon]$. Then $\alpha_1(1) = a_1$ and $\beta_1(1)
= b_1$.

Let $\phi: [0,1] \rightarrow [0,1]$, \ $ \phi \in C^1[0,1]$, be
such that
\begin{equation} \label{eqphi}
\phi(t) = 1, \ \ \ 0 \leq t \leq \epsilon/2, \ \ \ \ \ \ \phi(t) =
0, \ \ \ \epsilon \leq t \leq 1.
\end{equation}
Define the operator $S: L_{2,|r|}(\imath)\rightarrow
L_{2,|r|}(\jmath)$ by
\begin{equation} \label{eqdS}
(Sf)\bigl(\beta_1(t)\bigr) : = |\alpha'| f\bigl(\alpha_1(t)\bigr)
\phi(t), \ \ \ f \in L_{2,|r|}(\imath), \ \ \ t \in [0,1].
\end{equation}
Clearly $S$ is linear.

In what follows we shall use the combination of property
\eqref{eqphi} and Remark~\ref{rd0} to simplify the notation and
calculations. For example
these imply that in the definition \eqref{eqdS} of $S$ we could
use $\beta$ and $\alpha$ instead of $\beta_1$ and $\alpha_1$
without changing the substance of the definition.

At various points of the proof we shall employ the monotonic
(increasing or decreasing) substitutions
\[
x = \beta(t), \ \ \ \ \alpha(t) = \xi.
\]
To prove that $S$ is bounded we let $f \in
L_{2,|r|}(\imath)$ and calculate
\begin{align*}
\int_{\jmath} |(Sf)(x)|^2 |r(x)| dx
  & = \sgn(\beta') \int_{0}^{\epsilon}
  |(Sf)(\beta(t))|^2 |r(\beta(t))| \, \beta'  dt \\
& = |\alpha'|^2 |\beta'| \int_{0}^{\epsilon}
|f(\alpha(t))|^2\, |\phi(t)|^2 \, \rho(t) \,|r(\alpha(t))| dt  \\
& \leq |\alpha'|^2 |\beta'| R \int_{0}^{\epsilon}
|f(\alpha(t))|^2\, |\phi(t)|^2 \, |r(\alpha(t))| dx  \\
%
%
&  \leq |\alpha'|  |\beta'| R \int_{\imath}
|f(\xi)|^2\, |r(\xi)| d\xi,
\end{align*}
where $R$ is an upper bound of the function $\rho$.  The above
calculation proves that $S$ is bounded and $\|S\| \leq |\alpha'|
|\beta'| R$.

To verify the first claim in ($S$-\ref{itS2}), let $|x-b_1| < 1/2$
and $f \in L_{2,|r|}(\imath)$. Note that the length of $\jmath$ is
$1$, the endpoints of $\jmath$ are $b, b_1$, and $\beta(0) = b$.
Since $\beta\bigl([0,\epsilon]\bigr)$ has the length $ < 1/2$ and
since $\beta_1$ is strictly monotonic we conclude that $t =
\beta_1^{-1}(x) > \epsilon$.  Therefore, by \eqref{eqphi},
\[
(Sf)(x) = |\alpha'| f\bigl(\alpha(t)\bigr) \phi(t) = 0.
\]
This proves the first claim in ($S$-\ref{itS2}).

To prove $S \cF_{\max}(\imath) \subset  \cF_{\max}(\jmath)$, let
$f \in \cF_{\max}(\imath)$.  By definition \eqref{eqdS}, since $f$
is absolutely continuous on $\imath$ and $\phi \in C^1[0,1]$, the
function $Sf$ is absolutely continuous on $\jmath$ and for almost
all $t \in [0,1]$ we have
\begin{equation}
 \label{eqgSf'}
 \beta' (Sf)'(\beta(t))
 = |\alpha'| \bigl(\alpha'\, f'(\alpha(t)) \phi(t)
  + f(\alpha(t))\phi{'}(t) \bigr) .
\end{equation}
To prove that $Sf \in \cF_{\max}(\jmath)$ we need to show that
$(Sf)' \in L_{2,p}(\jmath)$, that is
\begin{equation}
  \label{eqgSf'i}
  \begin{split}
\int_{\jmath} |(Sf)'(x)|^2p(x)dx  =
  |\beta'| \int_0^\epsilon
  |(Sf)'(\beta(t))|^2 p(\beta(t))dt
 < + \infty.
  \end{split}
\end{equation}
We consider each summand in \eqref{eqgSf'} separately. By
\eqref{eqphi}, the second function in the sum in \eqref{eqgSf'} is
a continuous function which vanishes outside of the interval
$[\epsilon/2, \epsilon]$.  Since by assumption $p\circ\beta$ is an
integrable function on $[\epsilon/2, \epsilon]$, it follows that
\begin{equation}
 \label{eqgSf'st}
\int_0^1 |f(\alpha(t))\phi{'}(t)|^2 p(\beta(t)) dt < + \infty.
\end{equation}
Using the notation and assumptions from Definition~\ref{dscab},
for the first function in the sum in \eqref{eqgSf'} we have
\begin{align} \nonumber
\int_{0}^1 |f'(\alpha(t))|^2 |\phi(t)|^2 p(\beta(t)) dt & =
 \int_{0}^{\epsilon}
|f'(\alpha(t))|^2 \, |\phi(t)|^2 \, \varpi(t) \,p(\alpha(t)) dt
\\
\label{eqgSf'ft}
&  \leq \tau \int_{0}^{\epsilon}
 |f'(\alpha(t))|^2 \, |\phi(t)|^2 \,p(\alpha(t)) dt  \\ \nonumber
&  \leq  \frac{\tau}{|\alpha'|} \int_{\imath}
 |f'(\xi)|^2\,  p(\xi) d\xi .
\end{align}
Since $f' \in L_{2,p}(\imath)$ the last expression is finite.
Based on \eqref{eqgSf'}, \eqref{eqgSf'i}, \eqref{eqgSf'st} and
\eqref{eqgSf'ft} we conclude that $(S f)' \in L_{2,p}(\jmath)$ and
consequently $S f \in \cF_{\max}(\jmath)$.

The next step in the proof is to calculate
 $$
S^*: L_{2,|r|}(\jmath)\rightarrow L_{2,|r|}(\imath).
 $$
Note that $S^*$ is calculated with respect to the Hilbert space
inner products on the underlying spaces. Property \eqref{eqphi}
allows us to consider only affine changes of variable in the
integrals below. Let $f \in L_{2,|r|}(\imath)$ and $g \in
L_{2,|r|}(\jmath)$.  Then
\begin{align*}
\int_{\jmath} (Sf)(x) \, & \co{g(x)}\,|r(x)| dx \\
%
& = |\beta'| \int_{0}^{\epsilon}
  (Sf)(\beta(t)) \,\co{g(\beta(t))}\,|r(\beta(t))| \, dt \\
 & =  |\beta'| |\alpha'| \int_{0}^{\epsilon}
  f(\alpha(t))\, \phi(t)\,\co{g(\beta(t))}\,
   \rho(t) \,|r(\alpha(t))|\, dt  \\
& = |\beta'|  \int_{\imath}
 f(\xi)\, \phi(\alpha^{-1}(\xi))\,\co{g(\beta(\alpha^{-1}(\xi)))}\,
 \rho(\alpha^{-1}(\xi)) \,|r(\xi)|\, d\xi .
\end{align*}
Therefore for $g  \in L_{2,|r|}(\jmath)$ we have
\begin{equation*}
(S^*g)\bigl(x\bigr)
 : = |\beta'|
  \bigl( \rho \, \phi \, (g\circ \beta) \bigr)
    \bigl(\alpha^{-1}(x)\bigr),
 \ \ \ x \in \imath.
\end{equation*}
Thus
\begin{equation} \label{eqdS*}
(S^*g)\bigl(\alpha(t)\bigr)
  = |\beta'|  g\bigl(\beta(t)\bigr) \rho(t) \, \phi(t),
 \ \ \ g \in L_{2,|r|}(\jmath), \ \ \
t \in [0,1].
\end{equation}
As the adjoint of a bounded operator, the operator $S^*$
is bounded. To verify the second part of
($S$-\ref{itS2}) let $|x-a_1| < 1/2$ and $g \in
L_{2,|r|}(\jmath)$.  Note that the length of $\imath$ is
$1$ and $\alpha(0) = a$.  Since
$\alpha\bigl([0,\epsilon]\bigr)$ has length $ < 1/2$ and
since $\alpha_1$ is strictly monotonic we conclude that
$t = \alpha_1^{-1}(x) > \epsilon$.  Therefore, by
\eqref{eqphi},
\[
(S^*g)(x) = |\beta'| g\bigl(\beta(t)\bigr) \rho(t) \phi(t) = 0.
\]

To prove $S^* \cF_{\max}(\jmath) \subset
\cF_{\max}(\imath)$, let $g \in \cF_{\max}(\jmath)$.
Since $g$ and $\rho$ are absolutely continuous and $\phi
\in C^1[0,1]$, the function $(S^*g)\circ \alpha$ is
absolutely continuous on $[0,1]$.  Differentiation of
\eqref{eqdS*} yields
\begin{multline}
 \label{eqgS*f'}
\alpha' \, (S^*g)'\bigl(\alpha(t)\bigr) \\
  = |\beta'| \left(
  \bigl( \rho' \, \phi \, (g\circ \beta) \bigr)(t)
 + \bigl( \rho \, \phi' \, (g\circ \beta) \bigr)(t)
  + \beta' \bigl( \rho \, \phi \, (g'\circ \beta) \bigr) (t)
  \right),
\end{multline}
for almost all $t \in [0,1]$. To prove that $S^*f \in
\cF_{\max}(\imath)$ we need to show that $(S^*f)' \in
L_{2,p}(\imath)$, that is
\begin{equation}
  \label{eqgS*f'i}
  \begin{split}
\int_{\imath} |(S^*f)'(\xi)|^2p(\xi)d\xi
 & = |\alpha'| \int_0^\epsilon
  |(S^*f)'(\alpha(t))|^2 p(\alpha(t))dt < + \infty.
  \end{split}
\end{equation}
We prove that each summand on the right-hand side of
\eqref{eqgS*f'} belongs to $L_{2,p}(\imath)$.  By
\eqref{eqphi}, the second summand is a continuous
function which vanishes outside of the interval
$[\epsilon/2, \epsilon]$. Since $p\circ \alpha$ is an
integrable function on $[\epsilon/2, \epsilon]$, it
follows that
\begin{equation}
 \label{eqgS*f'st}
\int_0^1 |g(\beta(t))\phi{'}(t) \rho(t)|^2 p(\alpha(t)) dt < +
\infty.
\end{equation}
Next, we consider the third summand in \eqref{eqgS*f'}.
Since $\rho$ is continuous on $[0,1]$ we can consider
only $\phi \, (g'\circ \beta)$:
\begin{align} \nonumber
 \int_{0}^1
  \left|g'\bigl(\beta(t)\bigr) \,
  \phi(t)\right|^2 p(\alpha(t)) dt
  & = \int_{0}^{\epsilon}
  \left|g'\bigl(\beta(t)\bigr) \,
  \phi(t)\right|^2 \frac{1}{\varpi(t)} p(\beta(t)) dt
   \\ \label{eqgS*f'tt}
  & \leq  \tau   \int_{0}^{\epsilon}
  \left|g'\bigl(\beta(t)\bigr) \,
  \phi(t)\right|^2 p(\beta(t)) dt
 \\ \nonumber
&  = \frac{\tau}{ |\beta'|} \int_{\jmath} |g'(x)|^2  p(x) dx
< + \infty.
\end{align}
Finally, for the first summand in \eqref{eqgS*f'}, it is
sufficient to consider $\rho' \phi$, since $g\circ
\beta$ is absolutely continuous. By \eqref{eqdscabiii}
\begin{align} \label{eqgS*f'ft}
 \int_{0}^1
  \left|\rho'(t) \,
  \phi(t)\right|^2 p(\alpha(t)) dt
  & \leq \int_{0}^{\epsilon}
  \left|\rho'(t)\right|^2 p(\alpha(t)) dt < + \infty .
\end{align}
Based on \eqref{eqgS*f'}, \eqref{eqgS*f'i}, \eqref{eqgS*f'st},
\eqref{eqgS*f'tt} and \eqref{eqgS*f'ft} we conclude that $(S^* f)'
\in L_{2,p}(\imath)$ and consequently $S^* f \in
\cF_{\max}(\imath)$.

Thus we have verified the properties ($S$-\ref{itS1}),
($S$-\ref{itS2}), ($S$-\ref{itS4}). Since ($S$-\ref{itS3}) is
clear the theorem is proved.
\end{proof}

\section{Condition at $0$ and associated operator} \label{sc0}

\begin{condition}[Condition at $0$] \label{c0}
Let $p$ and $r$ be coefficients in \eqref{sl1}.  Denote by $0_{-}$
a generic left and by $0_{+}$ a generic right half-neighborhood of
$0$.  We assume that at least one of the four ordered pairs of
half-neighborhoods
\begin{equation} \label{eq4p}
(0_-,0_-), \ \ \ \ (0_-,0_+), \ \ \ \ (0_+,0_-), \ \ \ \
(0_+,0_+),
\end{equation}
is smoothly connected with connection parameters
$\alpha_0', \beta_0'$ and $\rho_0(0)$ such that
$|\alpha_0'| \neq |\beta_0'| \rho_0(0)$.
\end{condition}

\begin{theorem} \label{twat0}
Assume that the coefficients $p$ and $r$ satisfy Condition {\rm
\ref{c0}}.  Then there exists an operator $W_0:L_{2,r} \rightarrow
L_{2,r}$ such that
\begin{enumerate}[{\rm (a)}]
\item
\ $W_0$ is bounded on $L_{2,|r|}$.
\item \label{itwat0b}
The operator $J_0 W_0 - I$ is nonnegative on the Hilbert
space $L_{2,|r|}$.  In particular $W_0^{-1}$ is bounded
and $W_0$ is positive on the Krein space $L_{2,r}$.
\item \label{itwat0c}
$(W_0f)(x) = (Jf)(x), \ \ \ \frac{1}{2} \leq |x| \leq 1, \ \ \ f \in
L_{2,r}$.
\item
\ $W_0\cF_{\max}[-1,1] \subset \cF_{\max}[-1,1]$.
\end{enumerate}
\end{theorem}

\begin{proof}
Let $\alpha_0', \beta_0'$, and $\rho_0(0)$ be given by Condition
\ref{c0}. Recall that $|\alpha_0'| \neq |\beta_0'| \rho_0(0)$. Let
$\phi_0:[-1,1]\rightarrow [0,1],\, \phi_0 \in C^1[-1,1]$ be an
even function such that
\begin{equation*} 
\phi_0(0) = 1 \ \ \ \ \ \text{and} \ \ \ \ \ \phi_0(x) = 0 \ \ \
\text{for} \ \ \ 1/2 \leq |x| \leq 1.
\end{equation*}
Define the operators
\begin{equation*}
P_{0,-}: L_{2,|r|}(-1,0) \rightarrow L_{2,|r|}(-1,0) \ \ \
\text{and} \ \ \ P_{0,+}: L_{2,|r|}(0,1) \rightarrow
L_{2,|r|}(0,1)
\end{equation*}
by
\begin{alignat*}{2}
(P_{0,-} f)(x) &=  f(x)\, \phi_0(x), \ \ \ f \in L_{2,|r|}(-1,0),
 \ \ \  & & x \in [-1,0], \\
(P_{0,+} f)(x) &=  f(x)\, \phi_0(x), \ \ \ f \in L_{2,|r|}(0,1),
 & & x \in [0,1].
\end{alignat*}
Then $P_{0,-}$ and $P_{0,+}$ are self-adjoint operators
with the following properties:
\begin{alignat}{2} \label{eqpP1-}
(P_{0,-}f)(x) & = 0, \ \ \ f \in L_{2,|r|}(-1,0), \ \ \
 &  -1 \leq & \ x \leq -\tfrac{1}{2}, \\ \label{eqpP1+}
 (P_{0,+}f)(x) &= 0, \ \ \ f \in L_{2,|r|}(0,1), \ \ \
  & \tfrac{1}{2} \leq & \ x \leq 1,
\end{alignat}
\begin{equation} \label{eqpP2}
P_{0,-}\cF_{\max}[-1,0]  \subset \cF_{\max}[-1,0], \ \ \ \
 P_{0,+}\cF_{\max}[0,1]  \subset \cF_{\max}[0,1],
\end{equation}
and
\begin{align}
 \label{eqpP3-}
(P_{0,-}f)(0-) &= f(0-), \ \ \ \ f\in \cF_{\max}[-1,0], \\
  \label{eqpP3+}
(P_{0,+}f)(0+) & = f(0+), \ \ \ \ f\in  \cF_{\max}[0,1].
\end{align}
Here, the value of a function at $0\pm$ represents its one sided
limit.

Condition~\ref{c0} requires that one of the four ordered
pairs of half neighborhoods is smoothly connected. For
such a pair, Theorem~\ref{tgenS} guarantees the
existence of a specific operator which we denote by
$S_0$.  For each of the four pairs we shall use
different combinations of scaled operators $P_{0,-},
P_{0,+}, S_{0}$ and $S_{0}^*$ to define a bounded block
operator
\[
X_0 \ \ : \ \ \begin{matrix} L_{2,|r|}(-1,0)\\ \oplus \\
L_{2,|r|}(0,1)
\end{matrix} \ \ \rightarrow \ \
\begin{matrix} L_{2,|r|}(-1,0)\\ \oplus \\ L_{2,|r|}(0,1)
\end{matrix}
\]
with the following properties
\begin{align} \label{eqpX01}
(X_{0}^*f)(x) & = 0,  \ \ \ {1}/{2} \leq |x| \leq 1, \\
 \label{eqpX02}
X_{0}\cF_{\max}[-1,1] & \subset \cF_{\max}[-1,1], \\
 \label{eqpX02a}
(X_{0}f)(0) & = f(0) , \ \ \ f \in \cF_{\max}[-1,1], \\
 \label{eqpX03}
X_{0}^*\cF_{\max}[-1,1] & \subset \cF_{\max}[-1,0] \oplus
\cF_{\max}[0,1], \\
 \label{eqpX04}
 (X_{0}^*f)(0+)+(X_{0}^*f)(0-)& = - 2f(0) , \ \ \ f \in
\cF_{\max}[-1,1].
\end{align}
These properties of $X_0$ and $X_0^*$ imply that the operator
\[
W_0 = J\,\bigl(X_0^*X_0 + I\bigr)
\]
has all the properties stated in the theorem.

Since we assume that $|\alpha_0'| \neq |\beta_0'|
\rho_0(0)$, the system
\[
\gamma_1 |\alpha_0'| + \gamma_2 = 1, \ \ \ \ \ \
\gamma_1 |\beta_0'| \rho_0(0)  + \gamma_2 = -3
\]
has a nontrivial real solution $\gamma_1,\gamma_2$.  We
use this solution in the definitions below.

\noindent {\bf Case 1.} Assume that the half-neighborhoods
$0_-,0_-$ in \eqref{eq4p} are smoothly connected.  Then by Theorem
\ref{tgenS} there exists an operator
  $$
S_{0}: L_{2,|r|}(-1,0) \rightarrow L_{2,|r|}(-1,0)
  $$
which satisfies ($S$-\ref{itS1})-($S$-\ref{itS3}) in Theorem
\ref{tgenS} with $\imath = \jmath = [-1,0]$, $a = b = 0$. In
particular, for $f \in \cF_{\max}[-1,0]$,
\begin{equation*} 
(S_0f)(0-) = |\alpha_0'| \, f(0-), \ \ \ \ \ (S_0^*f)(0-) =
|\beta_0'|\, \rho_0(0) \, f(0-).
\end{equation*}
We define $X_0$ and calculate $X_0^*$ as
\begin{alignat*}{2}
X_0 & = \begin{bmatrix} \gamma_1 S_{0}  + \gamma_2 P_{0,-} & 0 \\
0 & P_{0,+}  \end{bmatrix}, & \ \ \ \ \ &
X_0^*  = \begin{bmatrix} \gamma_1 S_{0}^*  + \gamma_2 P_{0,-} & 0 \\
0 & P_{0,+}
\end{bmatrix}.
\end{alignat*}
{\bf Case 2.} Assume that the half-neighborhoods $0_-,0_+$ in
\eqref{eq4p} are smoothly connected. Then by Theorem~\ref{tgenS}
there exists an operator
 $$
S_{0}: L_{2,|r|}(-1,0) \rightarrow L_{2,|r|}(0,1)
 $$
which satisfies ($S$-\ref{itS1})-($S$-\ref{itS3}) in Theorem
\ref{tgenS} with $\imath = [-1,0],\, \jmath = [0,1]$, $a = b = 0$.
In particular, for $f \in \cF_{\max}[-1,0]$,
\begin{equation*} 
(S_0f)(0+) = |\alpha_0'| \, f(0-), \ \ \ \ \ (S_0^*f)(0-) =
|\beta_0'|\, \rho_0(0) \, f(0+).
\end{equation*}
We define $X_0$ and calculate $X_0^*$ as
\begin{alignat*}{2}
X_0 & = \begin{bmatrix} P_{0,-} & 0 \\ \gamma_1 S_{0} & \gamma_2
P_{0,+}
\end{bmatrix}, & \ \ \ \ \  &
 X_0^*  = \begin{bmatrix} P_{0,-} & \gamma_1 S_{0}^* \\ 0  &
 \gamma_2 P_{0,+}
\end{bmatrix}.
\end{alignat*}
{\bf Case 3.}  Assume that the half-neighborhoods $0_+,0_-$ in
\eqref{eq4p} are smoothly connected.  Then by Theorem~\ref{tgenS}
there exists an operator
 $$
S_{0}: L_{2,|r|}(0,1) \rightarrow L_{2,|r|}(-1,0)
 $$
which satisfies ($S$-\ref{itS1})-($S$-\ref{itS3}) in Theorem
\ref{tgenS} with $\imath = [0,1], \, \jmath = [-1,0]$, $a = b =
0$. In particular, for $f \in \cF_{\max}[0,1]$,
\begin{equation*} 
(S_0f)(0-) = |\alpha_0'| \, f(0+), \ \ \ \ (S_0^*f)(0+) =
|\beta_0'|\, \rho_0(0) \, f(0-).
\end{equation*}
We define $X_0$ and calculate $X_0^*$ as
\begin{alignat*}{2}
X_0 & = \begin{bmatrix} \gamma_2 P_{0,-} & \gamma_1 S_{0} \\
                              0 &  P_{0,+}
\end{bmatrix},& \ \ \ \ \ &
 X_0^*  = \begin{bmatrix} \gamma_2 P_{0,-} & 0 \\
                         \gamma_1 S_{0}^* &  P_{0,+}
\end{bmatrix}.
\end{alignat*}
{\bf Case 4.} Assume that the half-neighborhoods $0_+,0_+$ in
\eqref{eq4p} are smoothly connected.  Then by Theorem~\ref{tgenS}
there exists an operator
 $$
S_{0}: L_{2,|r|}(0,1) \rightarrow L_{2,|r|}(0,1)
 $$
which satisfies ($S$-\ref{itS1})-($S$-\ref{itS3}) in
Theorem~\ref{tgenS} with $\imath = [0,1], \, \jmath = [0,1]$, $a = b
= 0$. In particular, for $f \in \cF_{\max}[0,1]$,
\begin{equation*} 
(S_0f)(0+) = |\alpha_0'| \, f(0+), \ \ \ \ (S_0^*f)(0+)
 = |\beta_0'|\, \rho_0(0) \, f(0+).
\end{equation*}
We define $X_0$ and calculate $X_0^*$ as
\begin{alignat*}{2}
X_0 &= \begin{bmatrix} P_{0,-} & 0 \\
0 & \gamma_1 S_{0}+\gamma_2 P_{0,+}
\end{bmatrix}, & \ \ \ \ \  & X_0^*  = \begin{bmatrix}
P_{0,-} & 0 \\
0 & \gamma_1 S_{0}^* + \gamma_2 P_{0,+}
\end{bmatrix}.
\end{alignat*}

First note that in each of the four cases above, the operator
$X_0$ is bounded since each of its components is bounded.

In each of the four cases above, the property \eqref{eqpX01}
follows from ($S$-\ref{itS2}) in Theorem~\ref{tgenS}, and
properties \eqref{eqpP1-} and \eqref{eqpP1+}.

Let $f \in \cF_{\max}[-1,1]$.  Since $\gamma_1
|\alpha_0'| + \gamma_2 = 1$, the function $X_0f$ is
continuous in each case and $(X_0f)(0) = f(0)$.  This,
\eqref{eqpP2}, \eqref{eqpP3-}, \eqref{eqpP3+},
($S$-\ref{itS4}) and ($S$-\ref{itS3}) in Theorem
\ref{tgenS} imply \eqref{eqpX02}. Inclusion
\eqref{eqpX03} follows similarly.

In each of the above cases, equation \eqref{eqpX04} is a
consequence of
 $$
 \gamma_1 |\beta_0'| \rho_0(0) + \gamma_2 = -3,
 $$
\eqref{eqpP3-}, \eqref{eqpP3+} and ($S$-\ref{itS3}) in Theorem
\ref{tgenS}.  This proves the theorem.
\end{proof}

\begin{remark}
Note the behavior of the operator $W_0$ in Theorem~\ref{twat0} at
the boundary of the interval $[-1,1]$:
\[
\begin{bmatrix}
(W_0f)(-1) \\  (W_0f)(1)
\end{bmatrix} =
\begin{bmatrix}
- f(-1) \\  f(1)
\end{bmatrix}, \ \ \ f \in \cF_{\max}[-1,1].
\]
This property of $W_0$ will be used in
Section~\ref{sk=0}. In the next section, under
additional assumptions on the coefficients $p$ and $r$
in a neighborhood of $-1$ and $1$, we shall construct
operators $W$ with specified behaviors at $-1$ and $1$.
\end{remark}

\section{Conditions at -1 and 1, and associated operators}
\label{scm11}

In this section we show that under additional
assumptions on the coefficients $p$ and $r$ near $-1$ we
can construct an operator $W_{-1}$ with prescribed
behavior at $-1$ and under additional assumptions near
$1$ we can construct an operator $W_{+1}$ with
prescribed behavior at $1$.

\begin{condition}[Condition at $-1$] \label{cat-1}
Let $p$ and $r$ be coefficients in \eqref{sl1}.  We assume that a
right half-neighborhood of $-1$ is smoothly connected to a right
half-\!\! neighborhood of $-1$ with the connection parameters
$\alpha_{-1}', \beta_{-1}'$ and $\rho_{-1}(0)$ such that
$|\alpha_{-1}'| \neq |\beta_{-1}'| \rho_{-1}(0)$.
\end{condition}

\begin{condition}[Condition at $1$] \label{cat1}
Let $p$ and $r$ be coefficients in \eqref{sl1}.  We assume that a
left half-neighborhood of $1$ is smoothly connected to a left
half-neighborhood of $1$ with the connection parameters
$\alpha_{+1}', \beta_{+1}'$ and $\rho_{+1}(0)$ such that
$|\alpha_{+1}'| \neq |\beta_{+1}'| \rho_{+1}(0)$.
\end{condition}

In the rest of this section we shall need two operators
analogous to $P_{0,-}$ and $P_{0,+}$ introduced in
Section \ref{sc0}. Let $\phi_1:[-1,1] \rightarrow [0,1]$
be a smooth even function such that
\begin{equation} \label{eqphi1}
 \phi_1(-1) = 1, \ \ \ \ \ \ \
 \phi_1(x) = 0 \ \ \ \text{for} \ \ \  0 \leq |x| \leq 1/2,
 \ \ \ \ \ \ \  \phi_1(1) = 1.
\end{equation}
Define the operators
\begin{equation*}
P_{1,-}: L_{2,|r|}(-1,0) \rightarrow L_{2,|r|}(-1,0) \ \ \
\text{and} \ \ \  P_{1,+}: L_{2,|r|}(0,1) \rightarrow
L_{2,|r|}(0,1)
\end{equation*}
by
\begin{alignat}{2} \label{eqdP1-}
(P_{1,-} f)(x) &=  f(x)\, \phi_1(x), \ \ \ f \in L_{2,|r|}(-1,0),
\ \ \ & & x \in [-1,0], \\ \label{eqdP1+}
 (P_{1,+} f)(x) &= f(x)\, \phi_1(x), \ \ \ f \in L_{2,|r|}(0,1),
 \ \ \ \ & & x \in [0,1].
\end{alignat}
Then $P_{1,-}$ and $P_{1,+}$ are self-adjoint operators
with the following properties:
\begin{alignat}{2} \label{eqpP11-}
(P_{1,-}f)(x) & = 0, \ \ \ f \in L_{2,|r|}(-1,0), \ \ \
 &  -\tfrac{1}{2} \leq & \ x \leq 0, \\ \label{eqpP11+}
 (P_{1,+}f)(x) &= 0, \ \ \ f \in L_{2,|r|}(0,1), \ \ \
  & 0 \leq & \ x \leq \tfrac{1}{2},
\end{alignat}
\begin{equation} \label{eqpP12}
P_{1,-}\cF_{\max}[-1,0]  \subset \cF_{\max}[-1,0], \ \ \ \
 P_{1,+}\cF_{\max}[0,1]  \subset \cF_{\max}[0,1],
\end{equation}
and
\begin{align}
 \label{eqpP13-}
(P_{1,-}f)(-1+) &= f(-1+), \ \ \ \ f\in \cF_{\max}[-1,0], \\
  \label{eqpP13+}
(P_{1,+}f)(1-) & = f(1-), \ \ \ \ f\in  \cF_{\max}[0,1].
\end{align}

\begin{proposition} \label{pwat-1}
Assume that the coefficients $p$ and $r$ satisfy Condition {\rm
\ref{cat-1}}. Let $\mu$ be an arbitrary complex number.  Then
there exists an operator $W_{-1}:L_{2,r} \rightarrow L_{2,r}$ such
that
\begin{enumerate}[{\rm (a)\ }]
\item
 $W_{-1}$ is bounded on $L_{2,|r|}$.
\item \label{ipwat-1b}
The operator $J_0W_{-1} - I$ is nonnegative on the
Hilbert space $L_{2,|r|}$.  In particular
$(W_{-1})^{-1}$ is bounded and $W_{-1}$ is positive on
the Krein space $L_{2,r}$.
\item \label{ipwat-1c}
$(W_{-1}f)(x) = (Jf)(x), \ \ \ - \frac{1}{2} \leq x \leq
1, \ \ \ f \in L_{2,r}$.
\item
$W_{-1}\cF_{\max}[-1,1] \subset
\cF_{\max}[-1,0]\oplus\cF_{\max}[0,1]$.
\item
 $(W_{-1}f)(-1) = \mu f(-1) \ \ \ \text{for all} \ \ \ f \in
\cF_{\max}[-1,1]$.
\end{enumerate}
\end{proposition}

\begin{proof}
We use the notation introduced in Condition~\ref{cat-1}. By
Theorem~\ref{tgenS} there exists a bounded operator $S_{-1}:
L_{2,|r|}(-1,0) \rightarrow L_{2,|r|}(-1,0)$ such that
\begin{equation*}
S_{-1} \cF_{\max}[-1,0] \subset  \cF_{\max}[-1,0] \ \ \ \text{and}
\ \ \  S_{-1}^* \cF_{\max}[-1,0] \subset  \cF_{\max}[-1,0],
\end{equation*}
and, for all $f \in \cF_{\max}[-1,0]$,
 \begin{align*}
(S_{-1}f)(-1) & = |\alpha_{-1}'| \,f(-1), & (S_{-1}^*f)(-1) & =
|\beta_{-1}'|\, \rho_{-1}(0) f(-1).
\end{align*}
Let $\mu$ be an arbitrary complex number. Since we assume that
$|\alpha_{-1}'| \neq |\beta_{-1}'|\, \rho_{-1}(0)$, the complex
numbers $\gamma_1$ and $\gamma_2$ can be chosen such that
\[
\gamma_1 |\alpha_{-1}'| + \gamma_2 = 1, \ \ \ \ \ \ \co{\gamma}_1
|\beta_{-1}'|\, \rho_{-1}(0) + \co{\gamma}_2  = -\mu-1.
\]
We define $X_{-1}$ and calculate $X_{-1}^*$ as
\begin{align*}
X_{-1} & = \begin{bmatrix}
\gamma_1 S_{-1} + \gamma_2 P_{1,-} & 0 \\
0 & 0 \end{bmatrix}, &
 X_{-1}^* & =
\begin{bmatrix}
 \co{\gamma}_1 S_{-1}^* + \co{\gamma}_2 P_{1,-} & 0\\
0 & 0 \end{bmatrix}.
\end{align*}
Then for all $f \in \cF_{\max}[-1,1]$ we have
 \[
(X_{-1}f)(-1) = f(-1) \ \ \  \text{and} \ \ \  (X_{-1}^*f)(-1)=
(-\mu-1) f(-1).
 \]
Therefore
\[
W_{-1} = J \bigl(X_{-1}^*X_{-1}+I\bigr)
\]
has all the properties stated in the proposition.
\end{proof}

The proof of the next proposition is very similar to the preceding
proof, and will be omitted.

\begin{proposition} \label{pwat1}
Assume that the coefficients $p$ and $r$ satisfy Condition {\rm
\ref{cat1}}. Let $\mu$ be an arbitrary complex number.  Then there
exists an operator
 $$
 W_{+1}:L_{2,r} \rightarrow L_{2,r}
 $$
such that
\begin{enumerate}[{\rm (a)}]
\item
 $W_{+1}$ is bounded on $L_{2,|r|}$.
\item \label{ipwat1b}
The operator $J_0W_{+1} - I$ is nonnegative on the
Hilbert space $L_{2,|r|}$.  In particular
$(W_{+1})^{-1}$ is bounded and $W_{+1}$ is positive on
the Krein space $L_{2,r}$.
\item \label{ipwat1c}
$(W_{+1}f)(x) = (Jf)(x), \ \ \  - 1 \leq x \leq
\frac{1}{2}, \ \ \ f \in L_{2,r}$.
\item
 $W_{+1}\cF_{\max}[-1,1] \subset
\cF_{\max}[-1,0]\oplus\cF_{\max}[0,1]$.
\item
 $(W_{+1}f)(1) = \mu f(1) \ \ \ \text{for all} \ \ \ f \in
\cF_{\max}[-1,1]$.
\end{enumerate}
\end{proposition}

\section{Riesz basis of root vectors} \label{sk=0}

In this section we return to the eigenvalue problem
\eqref{sl1}-\eqref{bc1} and the operator $A$ associated with it.
We start with cases when the conditions in Section~\ref{scm11} are
not needed. We remark that the notation of Section~\ref{sop} is
used extensively in the rest of this section.

\begin{theorem} \label{tk0sebc}
Assume that the following three conditions are satisfied.
\begin{enumerate}[{\rm (a)}]
\item
The coefficients $p$ and $r$ satisfy Condition~{\rm~\ref{c0}}.
\item \label{itk0sebc2}
One of the following is true:
\begin{enumerate}[{\rm (i) \ }]
\item
 $\M{L}_n \neq \M{0}$,
\item
$ \M{L}= [1 \ \ 0 \ \ 0 \ \ 0]$,
\item
$\M{L}= [0 \ \ 1 \ \ 0 \ \ 0]$.
\end{enumerate}
\item \label{itk0sebc3}
One of the following is true:
\begin{enumerate}[{\rm (i) \ }]
\item
$\M{N}_n \neq \M{0}$,
\item
$\M{N}= [1 \ \ 0 \ \ 0 \ \ 0]$ \ and \ $\M{\Delta} < 0$,
\item
$\M{N}= [0 \ \ 1 \ \ 0 \ \ 0]$ \ and \ $\M{\Delta} > 0$.
\end{enumerate}
\end{enumerate}
Then there is a basis for each root subspace of $A$, so that the
union of all these bases is a Riesz basis of $L_{2,|r|} \oplus
\nC_{|\Delta|}$.
\end{theorem}

\begin{proof}
Assume first that $\M{N}_n \neq \M{0}$.  By \eqref{fd1}, the form
domain of $A$ when $\M{L}_n \neq \M{0}$ is given by
\begin{align*}
\fdom(A) & = \left\{ \begin{bmatrix} f \\ z
\end{bmatrix} \in
\begin{matrix}
L_{2,r} \\ \oplus \\ \nC_{\Delta} \end{matrix} \, : \,
f \in \cF_{\max}, \  z \in \nC \right\}, \\
\end{align*}
and in the other two cases in (\ref{itk0sebc2}),  \eqref{fd2}
gives
\[
\fdom(A) = \left\{ \begin{bmatrix} f \\ z
\end{bmatrix} \in
\begin{matrix}
L_{2,r} \\ \oplus \\ \nC_{\Delta} \end{matrix} \, : \, f \in
\cF_{\max}, \, z \in \nC, \,  \M{L}_e \V{b}_e(f) = 0 \right\},
\]
where $\M{L}_e \V{b}_e(f) = f(-1)$ in case (\ref{itk0sebc2}-ii) and
$\M{L}_e \V{b}_e(f) = f(1)$ in case (\ref{itk0sebc2}-iii).

Next assume (\ref{itk0sebc3}-ii).  Then $\M{N}_e \V{b}_e(f) =
f(-1)$ and \eqref{fd3} shows that
\begin{align*} 
\fdom(A) & = \left\{ \begin{bmatrix} f \\ f(-1)
\end{bmatrix} \in
\begin{matrix}
L_{2,r} \\ \oplus \\ \nC_{\Delta}
\end{matrix} \, : \, f \in \cF_{\max}\right\}
\end{align*}
when $\M{L}_n \neq \M{0}$, and in the other cases in
(\ref{itk0sebc2}), \eqref{fd4} gives
\[
 \fdom(A)  = \left\{ \begin{bmatrix} f \\ f(-1)
\end{bmatrix} \in
\begin{matrix}
L_{2,r} \\ \oplus \\ \nC_{\Delta} \end{matrix} \, : \, f
\in \cF_{\max}, \  \M{L}_e \V{b}_e(f) = 0\right\}.
\]

Finally, assume (\ref{itk0sebc3}-iii).  Then $\M{N}_e \V{b}_e(f) =
f(1)$ and \eqref{fd3} shows that
\begin{align*} 
\fdom(A) & = \left\{ \begin{bmatrix} f \\ f(1)
\end{bmatrix} \in
\begin{matrix}
L_{2,r} \\ \oplus \\ \nC_{\Delta}
\end{matrix} \, : \, f \in \cF_{\max}\right\}
\end{align*}
when $\M{L}_n \neq \M{0}$, and in the other cases in
(\ref{itk0sebc2}), \eqref{fd4} gives
\[
 \fdom(A)  = \left\{ \begin{bmatrix} f \\ f(1)
\end{bmatrix} \in
\begin{matrix}
L_{2,r} \\ \oplus \\ \nC_{\Delta} \end{matrix} \, : \, f \in
\cF_{\max}, \  \M{L}_e \V{b}_e(f) = 0\right\}.
\]

Let $W_0$ be the operator constructed in Theorem~\ref{twat0}, and
let
\[
W = \begin{bmatrix} W_0 \ & 0 \\  0 \ & \sgn(\M{\Delta})
\end{bmatrix} \ : \
 \begin{matrix}L_{2,r} \\ \oplus  \\ \nC_{\M{\Delta}}
\end{matrix}
\ \ \rightarrow \ \ \begin{matrix}L_{2,r} \\ \oplus  \\
\nC_{\M{\Delta}} \end{matrix}.
\]
A straightforward verification shows that $W$ is a bounded,
boundedly invertible, positive operator in the Krein space
$L_{2,r} \oplus \nC_{\Delta}$ and $W\fdom(A) \subset \fdom(A)$ in
each of the above listed cases. Consequently, the theorem follows
from Theorem~\ref{W}.
\end{proof}

In the next result we shall assume that one of the conditions from
Section~\ref{scm11} is satisfied.

\begin{theorem} \label{td0k04}
Assume that the following three conditions are satisfied.
\begin{enumerate}[{\rm (a)}]
\item
The coefficients $p$ and $r$ satisfy Condition~{\rm~\ref{c0}}, and
one of Conditions~{\rm~\ref{cat-1}}, {\rm~\ref{cat1}}.
\item
  $\M{N}_n \neq \M{0}$.
\item
$\M{L} = [u \ \ v \ \ 0 \ \ 0]$ with $uv \neq 0$.
\end{enumerate}
Then there is a basis for each root subspace of $A$, so that the
union of all these bases is a Riesz basis of $L_{2,|r|} \oplus
\nC_{|\Delta|}$.
\end{theorem}

\begin{proof}
Under the assumptions of the theorem, \eqref{fd2} shows that the
form domain of $A$ is given by
\begin{equation*}
\fdom(A)  = \left\{ \begin{bmatrix} f \\ z
\end{bmatrix} \in
\begin{matrix}
L_{2,r} \\ \oplus \\ \nC_{\Delta} \end{matrix} \, : \, f \in
\cF_{\max}, \, z \in \nC, \, uf(-1) + vf(1) = 0 \right\}.
\end{equation*}
Define the following two Krein spaces:
\begin{equation} \label{eqK0K1}
 \cK_0  := L_{2,r}\!\left(-\tfrac{1}{2},\tfrac{1}{2}\right),  \ \ \
 \ \ \   \cK_1 := L_{2,r}(-1,-\tfrac{1}{2})[\dot{+}]
 L_{2,r}(\tfrac{1}{2},1).
\end{equation}
Extending the functions in $\cK_0$ and $\cK_1$ by $0$
onto the rest of $[-1,1]$, we can consider  the spaces
$\cK_0$ and $\cK_1$ as subspaces of $L_{2,r}$.

Then
 \[
L_{2,r} = \cK_0[\dot{+}]\cK_1.
 \]

Assume that the functions $p$ and $r$ satisfy
Conditions~\ref{c0} and \ref{cat-1}. Let $W_0$ be the
operator constructed in Theorem \ref{twat0} and let
$W_{-1}$ be the operator constructed in
Proposition~\ref{pwat-1} with $\mu = 1$.  Then
properties (\ref{ipwat-1c}) in Theorem~\ref{twat0} and
Proposition~\ref{pwat-1}, imply that $\cK_0$ and $\cK_1$
are invariant under $W_0$ and $W_{-1}$.  As we chose
$\mu = 1$, we have $(W_{-1}f)(-1) = f(-1)$ and
$(W_{-1}f)(1) = f(1)$. Define
\begin{equation} \label{eqdW011}
W_{01} : = W_0|_{_{\cK_0}}[\dot{+}]W_{-1}|_{_{\cK_1}}.
\end{equation}
Since $W_0$ and $W_{-1}$ are bounded, boundedly invertible and
positive in the Krein space $L_{2,r}$, so is the the operator
$W_{01}$.  Also, $W_{01} \cF_{\max}[-1,1] \subset
\cF_{\max}[-1,1]$ and
\begin{equation} \label{eqW01-}
\begin{bmatrix} (W_{01}f)(-1) \\ (W_{01}f)(1)
\end{bmatrix} = \begin{bmatrix} f(-1)
\\ f(1) \end{bmatrix}.
\end{equation}

If the functions $p$ and $r$ satisfy Conditions~\ref{c0} and
\ref{cat1}, then, instead of $W_{-1}$, we use the operator
$W_{+1}$ constructed in Proposition~\ref{pwat1} with $\mu = -1$.
Redefining the operator $W_{01}$ as
\begin{equation} \label{eqdW012}
W_{01} : = W_0|_{_{\cK_0}}[\dot{+}]W_{+1}|_{_{\cK_1}}
\end{equation}
we see that it is again bounded, boundedly invertible, and positive
in the Krein space $L_{2,r}, \; W_{01} \cF_{\max}[-1,1] \subset
\cF_{\max}[-1,1]$ and (since we use $\mu = -1$)
\begin{equation} \label{eqW01+}
\begin{bmatrix} (W_{01}f)(-1) \\ (W_{01}f)(1)
\end{bmatrix} = - \begin{bmatrix} f(-1)
\\ f(1) \end{bmatrix}.
\end{equation}

Now a simple inspection shows that, in both above cases, the
operator
\[
W = \begin{bmatrix} W_{01} \ & 0 \\  0 \ &
\M{\Delta}
\end{bmatrix} \ : \
 \begin{matrix}L_{2,r} \\ \oplus  \\ \nC_{\M{\Delta}}
\end{matrix}
\ \ \rightarrow \ \ \begin{matrix}L_{2,r} \\ \oplus  \\
\nC_{\M{\Delta}} \end{matrix}
\]
is bounded, boundedly invertible and positive in the Krein space
$L_{2,r} \oplus \nC_{\Delta}$ and $W\fdom(A) \subset \fdom(A)$.
Thus the theorem again follows from Theorem~\ref{W}.
\end{proof}

Our final result covers the remaining cases, but in the interests
of presentation we shall impose no conditions on $\M{L}$. Of
course, there is some overlap with Theorem~\ref{tk0sebc}.

\begin{theorem} \label{tn=k2}
Assume that the following two conditions are satisfied.
\begin{enumerate}[{\rm (a)}]
\item  \label{tn=k22}
The coefficients $p$ and $r$ satisfy Condition~{\rm~\ref{c0}}, and
\begin{enumerate}[{\rm (i)}]
\item
Condition~{\rm~\ref{cat-1}} if $\M{\Delta} > 0$,
\item
Condition~{\rm~\ref{cat1}} if $\M{\Delta} < 0$.
\end{enumerate}
\item \label{tn=k21}
$\M{N}_n = \M{0}$.
\end{enumerate}
Then there is a basis for each root subspace of $A$, so that the
union of all these bases is a Riesz basis of $L_{2,|r|} \oplus
\nC_{|\Delta|}$.
\end{theorem}

\begin{proof}
In this case \eqref{fd3} shows that the form domain of $A$ is
\begin{equation*}
 \fdom(A)  = \left\{ \begin{bmatrix} f \\ \M{N}_e \V{b}_e(f)
\end{bmatrix} \in
\begin{matrix}
L_{2,r} \\ \oplus \\ \nC_{\Delta}
\end{matrix} \, : \, f \in \cF_{\max} \right\}
\end{equation*}
if $\M{L}_n \ne \M{0}$, and
\begin{equation*}
 \fdom(A)  = \left\{ \begin{bmatrix} f \\ \M{N}_e \V{b}_e(f)
\end{bmatrix} \in
\begin{matrix}
L_{2,r} \\ \oplus \\ \nC_{\Delta}
\end{matrix} \, : \, f \in \cF_{\max}, \
  \M{L}_e \V{b}_e(f) = \M{0}\right\}
\end{equation*}
by \eqref{fd4} if $\M{L}_n = \M{0}$.

Let $W_{01}:L_{2,r}\to L_{2,r}$ be the operator constructed in the
proof of Theorem~\ref{td0k04}, and define $W: L_{2,r} \oplus
\nC_{\M{\Delta}} \to L_{2,r} \oplus \nC_{\M{\Delta}}$ by
\[
W = \begin{bmatrix} W_{01} \ & 0 \\  0 \ & \sgn(\M{\Delta})
\end{bmatrix}.
\]
As before, in both cases the properties of $W_{01}$ imply that $W$
is a bounded, boundedly invertible, positive operator in the Krein
space $L_{2,r} \oplus \nC_{\Delta}$ and $W\fdom(A) \subset
\fdom(A)$. Now the theorem follows from Theorem~\ref{W}.
\end{proof}

\begin{remark}
It is instructive to look at the above results from the viewpoint of
non-essential boundary conditions ($\M{L}_n,\,\M{N}_n \ne \M{0}$)
and essential ones (whose essential parts $\M{L}_e,\,\M{N}_e$ can be
separated or not). Let us call a boundary condition essentially
separated if it is either non-essential, or else its essential part
is separated. Theorem~\ref{tk0sebc} states that if both boundary
conditions are essentially separated, then subject to the sign
conditions in (\ref{itk0sebc3}-ii) and (\ref{itk0sebc3}-iii),
Condition~\ref{c0} suffices for the existence of a Riesz basis of
root vectors.

If any of these assumptions fail, then we impose conditions from
Section~\ref{scm11}. In particular, if the $\lambda$-dependent
boundary condition is non-essential, then either of these
conditions suffice, but in other cases the choice is governed by
the sign of $\Delta$.
\end{remark}

We conclude with a simple example.
\begin{example}
We suppose that
\[ \M{L} = [d_1 \ \ d_2 \ \ d_3 \ \ d_4],\;\;
   \M{M} = [m_1 \ \ m_2 \ \ m_3 \ \ m_4],\;\;
   \M{N} = [0 \ \ \gamma \ \ 0 \ \ 0],
\]
where $(d_3, \, d_4) \ne (0,0)$ and $\gamma m_4 > 0$.
Note that the only $\lambda$-dependent term in \eqref{bc1}
involves $f(1)$.

We calculate
\[ \M{M} \M{Q}^{-1} \M{N}^* = -i \gamma m_4, \]
so by \eqref{Del}, $\Delta > 0$. It follows from Theorem
\ref{tk0sebc} with condition (b-i) that Condition~\ref{c0} suffices
for the existence of a Riesz basis of root vectors.

This example overlaps with \cite[Corollary~3.8]{F3}, where
separated boundary conditions, also satisfying $d_2 = d_4 = m_1 =
m_3 = 0,\, d_3 = m_4 = 1,\, \gamma > 0$, were considered by Fleige
for a Krein-Feller equation instead of \eqref{sl1}.
\end{example}

\subsection*{Acknowledgment}
We thank a referee for a very careful reading of the submitted
version of this article. This has led to the correction of a number
of inaccuracies.

\end{document}